\documentclass[a4paper,12pt,authoryear]{elsarticle} 
\usepackage[left=3.0cm,right=2.0cm,top=2cm,bottom=3cm,bindingoffset=0cm]{geometry}

\usepackage{graphicx}
\usepackage{natbib}
\usepackage[english]{babel}

\usepackage{tabularx}
\usepackage{multirow}
\usepackage{multicol}
\usepackage{amssymb}
\usepackage{amstext}
\usepackage{amsmath}
\usepackage{pslatex}
\usepackage{floatrow}

\usepackage[update,prepend]{epstopdf}

\newtheorem{definition}{Definition}

\usepackage{lmodern}

\journal{arXiv}

\begin{document}

\begin{frontmatter}

\title{
      Lock-in range of PLL-based circuits with \\
      proportionally-integrating filter and \\ sinusoidal phase detector characteristic
}

\author{K.~D. Aleksandrov}
\author{N.V. Kuznetsov\corref{cor}\,}
\ead{nkuznetsov239@gmail.com}
\author{G.~A. Leonov}
\author{M.~V. Yuldashev}
\author{R.~V. Yuldashev}

\address{Faculty of Mathematics and Mechanics,\
Saint-Petersburg State University, Russia}
\address{Dept. of Mathematical Information Technology,\
University of Jyv\"{a}skyl\"{a}, Finland}
\address{Institute of Problems of Mechanical Engineering RAS, Russia}

\begin{abstract}
In the present work PLL-based circuits with sinusoidal phase detector characteristic and active proportionally-integrating (PI) filter are considered. The notion of lock-in range -- an important characteristic of PLL-based circuits, which corresponds to the synchronization without cycle slipping, is studied. For the lock-in range a rigorous mathematical definition is discussed. Numerical and analytical estimates for the lock-in range are obtained.
\end{abstract}

\end{frontmatter}




\section{Model of PLL-based circuits in the signal's phase space}
\label{sec:PLLBased}
For the description of PLL-based circuits a physical model in the signals space and a mathematical model in the signal's phase space are used \citep{Gardner-1966,ShahgildyanL-1966,Viterbi-1966}.

The equations describing the model of PLL-based circuits in the signals space are difficult for the study, since that equations are nonautonomous (see, e.g., \citep{KudrewiczW-2007}). By contrast, the equations of model in the signal's phase space are autonomous
\citep{Gardner-1966,ShahgildyanL-1966,Viterbi-1966}, what simplifies the study of PLL-based circuits.
The application of averaging methods \citep{MitropolskyB-1961,Samoilenko-2004-averiging} allows one to reduce the model of PLL-based circuits in the signals space to the model in
the signal's phase space
(see, e.g., \citep{LeonovKYY-2012-TCASII,LeonovK-2014-book,LeonovKYY-2015-SIGPRO,KuznetsovLSYY-2015-PD,KuznetsovKLNYY-2015-ISCAS,BestKKLYY-2015-ACC}.
\begin{figure}[!htbp]
\centering
\includegraphics[width=0.9\textwidth]{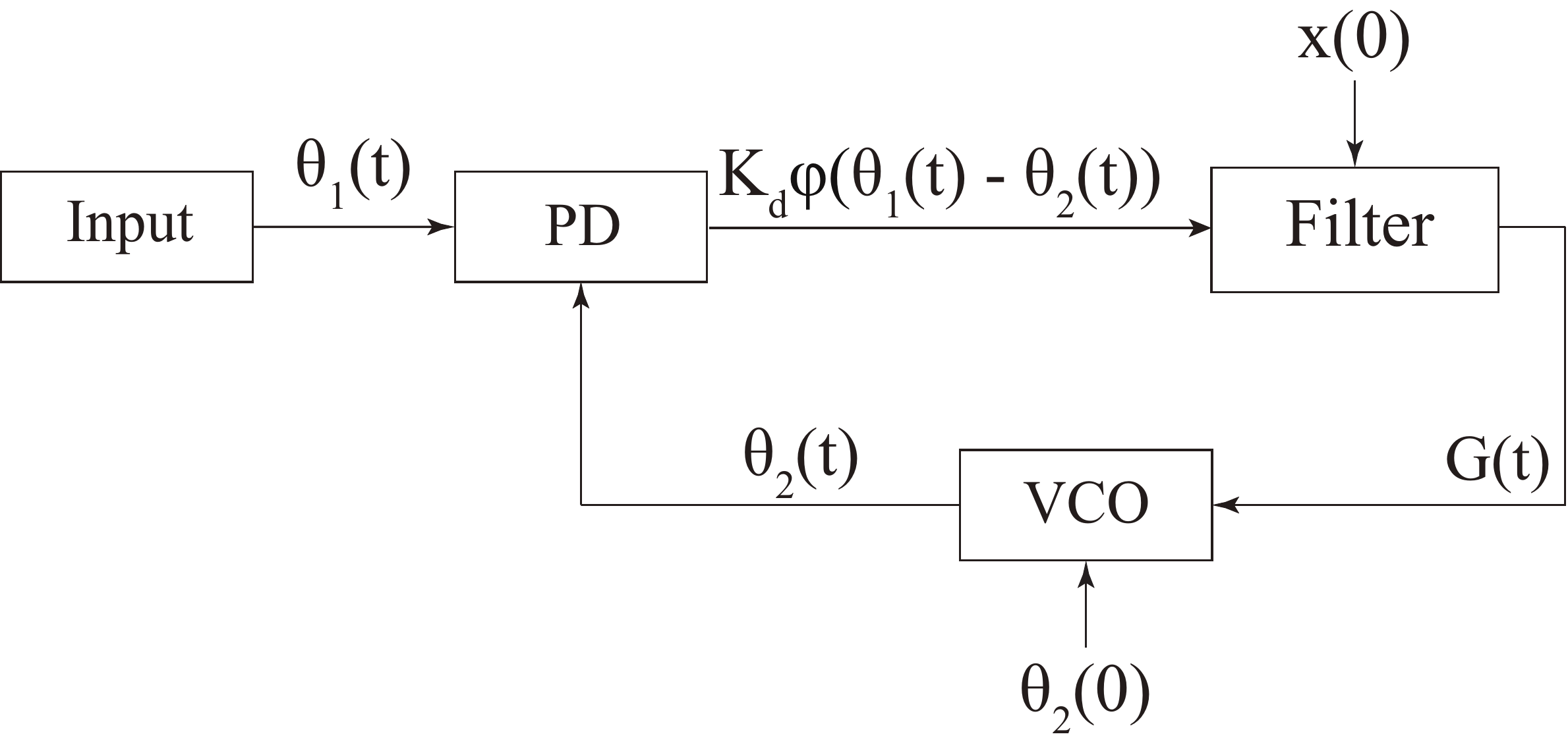}
\caption{Model of PLL-based circuit in the signal's phase space.}
\label{ris:PLLBased}
\end{figure}

Consider a model of PLL-based circuits in the signal's phase space (see Fig.~\ref{ris:PLLBased}).
A reference oscillator (Input) and a voltage-controlled oscillator (VCO) generate phases
$\theta_1(t)$ and $\theta_2(t)$, respectively.
The frequency of reference signal usually assumed to be constant:
\begin{equation}
\dot{\theta}_1(t) = \omega_1.
\label{eq:Input}
\end{equation}
The phases $\theta_1(t)$ and $\theta_2(t)$ enter the inputs of the phase detector (PD). The output of the phase detector in the signal's phase space is called a phase detector characteristic and has the form $$ K_d \varphi(\theta_1(t) - \theta_2(t)).$$ The maximum absolute value of PD output $K_d > 0$ is called a phase detector gain (see, e.g., \citep{Best-2007, Goldman-2007-book}). The periodic function $\varphi(\theta_\Delta(t))$ depends on difference $\theta_1(t) - \theta_2(t)$ (which is called a phase error and denoted by $\theta_\Delta(t)$). The PD characteristic depends on the design of PLL-based circuit and the signal waveforms $f_1(\theta_1)$ of Input and $f_2(\theta_2)$ of VCO.
In the present work a sinusoidal PD characteristic with $$\varphi(\theta_\Delta(t)) = \sin(\theta_\Delta(t))$$ is considered (which corresponds, e.g., to the classical PLL with $f_1(\theta_1(t)) = \sin(\theta_1(t))$ and $f_2(\theta_2(t)) = \cos(\theta_2(t))$).

The output of phase detector is processed by Filter. Further we consider the active PI filter (see, e.g., \citep{Baker-2011}) with transfer function
$W(s) = \frac{1 + \tau_2 s}{\tau_1 s}$, $\tau_1~>~0$, $\tau_2~>~0$. The considered filter can be described as
\begin{equation}
\begin{cases}
\dot{x}(t) = K_d \sin(\theta_\Delta(t)), \\
G(t) = \frac{1}{\tau_1}x(t) + \frac{\tau_2}{\tau_1}K_d\sin(\theta_\Delta(t)),
\label{eq:FilterOrig}
\end{cases}
\end{equation}
where $x(t)$ is the filter state.
The output of Filter $G(t)$ is used as a control signal for VCO:
\begin{equation}
\dot{\theta}_2(t) = \omega_2^{\rm free} + K_{v} G(t),
\label{eq:VCO}
\end{equation}
where $\omega_2^{\rm free}$ is the VCO free-running frequency and $K_{v} > 0$ is the VCO gain coefficient.

Relations (\ref{eq:Input}), (\ref{eq:FilterOrig}), and (\ref{eq:VCO}) result in autonomous system of differential equations
\begin{equation}
\begin{cases}
\dot{x} = K_d \sin(\theta_\Delta), \\
\dot{\theta}_\Delta = \omega_1 - \omega_2^{\rm free} - \frac{K_{v}}{\tau_1}\left(x + \tau_2 K_d \sin(\theta_\Delta)\right).
\end{cases}
\label{sys:PLLSys_beforeChange}
\end{equation}
Denote the difference of the reference frequency and the VCO free-running frequency $\omega_1 - \omega_2^{\rm free}$ by $\omega_\Delta^{\rm free}$.
By the linear transformation $x \rightarrow K_d x$ we have
\begin{equation}
\begin{cases}
\dot{x} = \sin(\theta_\Delta), \\
\dot{\theta}_\Delta = \omega_\Delta^{\rm free} - \frac{K_0}{\tau_1}\left(x + \tau_2 \sin(\theta_\Delta)\right),
\end{cases}
\label{sys:PLLSys}
\end{equation}
where $K_0 = K_{v} K_d$ is the loop gain. For signal waveforms listed in Table \ref{table:PDgains}, relations (\ref{sys:PLLSys}) describe the models of the classical PLL and two-phase PLL in the signal's phase space. The models of classical Costas loop and two-phase Costas loop in the signal's phase space can be described by relations similar to (\ref{sys:PLLSys}) (PD characteristic of the circuits usually is a $\pi$-periodic function, and the approaches presented in this paper can be applied to these circuits as well)
(see, e.g., \citep{BestKLYY-2014-IJAC,LeonovKYY-2015-SIGPRO,BestKKLYY-2015-ACC}). %
\begin{table}[ht]
\begin{center}
\renewcommand{\arraystretch}{1.5}
\begin{tabular}{|l|c|}
\hline
Signal waveforms & PD characteristic\\
\hline
\multicolumn{2}{|c|}{Classical PLL}\\
\hline
$f_{1}(\theta_1) = \sin(\theta_1)$&\multirow{2}{*}{$\frac{1}{2} \sin(\theta_\Delta)$}\\
$f_{2}(\theta_2) = \cos(\theta_2)$&\\
\hline
$f_{1}(\theta_1) = \sin(\theta_1)$&\multirow{2}{*}{$\frac{2}{\pi} \sin(\theta_\Delta)$}\\
$f_{2}(\theta_2) = \operatorname{sign}(\cos(\theta_2))$&\\
\hline
$f_{1}(\theta_1) = \begin{cases} \frac{2}{\pi}\theta_1 + 1, \theta_1 \in \left[0; \pi \right], \\ 1-\frac{2}{\pi}\theta_1, \theta_1 \in \left[\pi; 2\pi \right]
\end{cases}$&\multirow{2}{*}{$\frac{4}{\pi^2} \sin(\theta_\Delta)$} \\
$f_{2}(\theta_2) = \sin(\theta_2)$& \\
\hline
\multicolumn{2}{|c|}{Two-phase PLL}\\
\hline
$f_{1}(\theta_1) = \cos(\theta_1)$&\multirow{2}{*}{$\sin(\theta_\Delta)$}\\
$f_{2}(\theta_2) = \cos(\theta_2)$&\\
\hline
\end{tabular}
\end{center}
\caption{The dependency PD characteristics of PLL-based circuits on signal waveforms.}
\label{table:PDgains}
\end{table}

By the transformation $$\left(\omega_\Delta^{\rm free}, x, \theta_\Delta\right) \rightarrow \left(-\omega_\Delta^{\rm free}, -x, -\theta_\Delta\right),$$ (\ref{sys:PLLSys}) is not changed.
This property allows one to use the concept of frequency deviation $$\left|\omega_\Delta^{\rm free}\right| = \left|\omega_1 - \omega_2^{\rm free}\right|$$ and consider (\ref{sys:PLLSys}) with $\omega_\Delta^{\rm free} > 0$ only.

The state of PLL-based circuits for which the VCO frequency is adjusted to the reference frequency of Input is called a locked state.
The locked states correspond to the locally asymptotically stable equilibria of (\ref{sys:PLLSys}), which can be found from the relations
\begin{equation}
\begin{cases}
\sin(\theta_{eq}) = 0, \\
\omega_\Delta^{\rm free} - \frac{K_0}{\tau_1} x_{eq} = 0.
\end{cases}
\nonumber
\label{sys:PLLSysEq}
\end{equation}
Here $x_{eq}$ depends on $\omega_\Delta^{\rm free}$ and further is denoted by $x_{eq}(\omega_\Delta^{\rm free})$.

Since (\ref{sys:PLLSys}) is $2\pi$-periodic in $\theta_\Delta$, we can consider (\ref{sys:PLLSys}) in a $2\pi$-interval of $\theta_\Delta$, $\theta_\Delta \in \left(-\pi, \pi\right]$.
In interval $\theta_\Delta \in \left(-\pi, \pi\right]$ there exist two equilibria: $\left(0, \frac{\omega_\Delta^{\rm free} \tau_1}{K_0}\right)$ and $\left(\pi, \frac{\omega_\Delta^{\rm free} \tau_1}{K_0}\right)$. To define type of the equilibria let us write out corresponding characteristic polynomials and find the eigenvalues:
\begin{equation}
\begin{array}{ll}
\text{equilibrium }(0, \frac{\omega_\Delta^{\rm free} \tau_1}{K_0}): & \lambda^2 + \frac{K_0 \tau_2}{\tau_1}\lambda + \frac{K_0}{\tau_1} = 0; \\
\\
\lambda_{1,2} = \frac{-K_0 \tau_2 \pm \sqrt{(K_0 \tau_2)^2 - 4 K_0 \tau_1}}{2 \tau_1}, & (K_0 \tau_2)^2 - 4 K_0 \tau_1>0; \\
\lambda_{1} = \lambda_{2} = \frac{-K_0 \tau_2}{2 \tau_1}, \hskip0.2cm & (K_0 \tau_2)^2 - 4 K_0 \tau_1=0;\\
\lambda_{1,2} = \frac{-K_0 \tau_2 \pm i \sqrt{4 K_0 \tau_1 - (K_0 \tau_2)^2}}{2 \tau_1}, & (K_0 \tau_2)^2 - 4 K_0 \tau_1<0;\\
\end{array} \nonumber
\end{equation}
\\
\begin{equation}
\begin{array}{ll}
\text{equilibrium }(\pi, \frac{\omega_\Delta^{\rm free} \tau_1}{K_0}): & \lambda^2 - \frac{K_0 \tau_2}{\tau_1}\lambda - \frac{K_0}{\tau_1} = 0;\\
\\
\lambda_{1,2} = \frac{K_0 \tau_2 \pm \sqrt{(K_0 \tau_2)^2 + 4 K_0 \tau_1}}{2 \tau_1}. &
\end{array} \nonumber
\end{equation}

Denote the stable equilibrium as
\[
  \left(\theta_{eq}^s, x_{eq}(\omega_\Delta^{\rm free})\right) = \left(0, \frac{\omega_\Delta^{\rm free} \tau_1}{K_0}\right)
\]
and the unstable equilibrium as
\[
\left(\theta_{eq}^u, x_{eq}(\omega_\Delta^{\rm free})\right) = \left(\pi, \frac{\omega_\Delta^{\rm free} \tau_1}{K_0}\right).
\]
Thus, for any arbitrary $\omega_\Delta^{\rm free}$ the equilibria $$\left(\theta_{eq}^s + 2\pi k, x_{eq}(\omega_\Delta^{\rm free})\right) = \left(2\pi k, \frac{\omega_\Delta^{\rm free} \tau_1}{K_0}\right)$$ are locally asymptotically stable. Hence, the locked states of (\ref{sys:PLLSys}) are given by equilibria $\left(\theta_{eq}^s + 2 \pi k, x_{eq}(\omega_\Delta^{\rm free})\right)$. The remaining equilibria $$\left(\theta_{eq}^u + 2\pi k, x_{eq}(\omega_\Delta^{\rm free})\right) = \left(\pi + 2\pi k, \frac{\omega_\Delta^{\rm free} \tau_1}{K_0}\right)$$ are unstable saddle equilibria.

\section{The global stability of PLL-based circuit model}
\label{sec:StabilitySets}
In order to consider the lock-in range of PLL-based circuits let us discuss the global asymptotic stability. \textit{If for a certain $\omega_\Delta^{\rm free}$ any solution of (\ref{sys:PLLSys}) tends to an equilibrium, then the system with such $\omega_\Delta^{\rm free}$ is called globally asymptotically stable} (see, e.g., \citep{LeonovKYY-2015-TCAS}).
To prove the global asymptotic stability of (\ref{sys:PLLSys}) two approaches can be applied: the phase plane analysis \citep{Tricomi-1933,AndronovVKh-1937} and construction of the Lyapunov functions \citep{Lyapunov-1892}.

\begin{figure*}[!htbp]
\centering
\includegraphics[width=\linewidth]{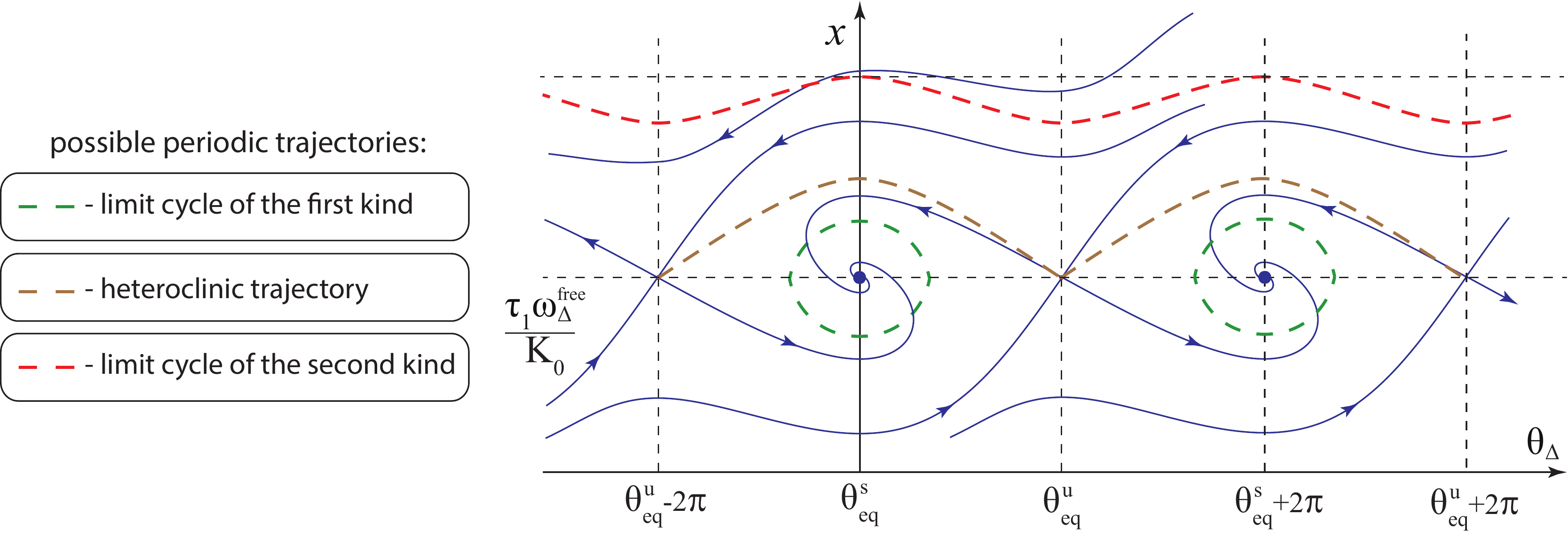}
\caption{Phase portrait and possible periodic trajectories of (\ref{sys:PLLSys}).}
\label{ris:CyclesPhasePlane}
\end{figure*}

By methods of the phase plane analysis, in \citep{Viterbi-1966} the global asymptotic stability of (\ref{sys:PLLSys}) for any $\omega_\Delta^{\rm free}$ is stated.
However, to complete rigorously the proof given in \citep{Viterbi-1966}, the additional explanations are required (i.e., the absence of heteroclinic trajectory and limit cycles of the first kind (see Fig.~\ref{ris:CyclesPhasePlane}) is needed to be explained; e.g., for the case of lead-lag filter a number of works \citep{Kapranov-1956,Gubar-1961,Shakhtarin-1969,Belyustina-1970} is devoted to the study of these periodic trajectories).

To overcome these difficulties, the methods of the Lyapunov functions construction can be applied. The modifications of the classical global stability criteria for cylindrical phase space are developed in \citep{GeligLY-1978,LeonovK-2014-book,LeonovKYY-2015-TCAS}.
The global asymptotic stability of (\ref{sys:PLLSys}) for any $\omega_\Delta^{\rm free}$
can be using the Lyapunov function
\begin{align*}
&V(x, \theta_\Delta) = \frac{1}{2} \left(x - \frac{\tau_1 \omega_\Delta^{\rm free}}{K_0}\right)^2 + \frac{2\tau_1}{K_0} \sin^2 \left(\frac{\theta_\Delta}{2}\right) \geq 0; \\
& \dot{V}(x, \theta_\Delta) = - \tau_2 \sin^2 \left(\theta_\Delta \right) < 0, \quad \forall \theta_\Delta \neq \{\theta^s_{eq} + 2 \pi k,\theta^u_{eq} + 2 \pi k\}.
\end{align*}



\section{The lock-in range definition and analysis}
\label{sec:LockIn}
Since the considered model of PLL-based circuits in the signal's phase space is globally asymptotically stable, it achieves locked state for any initial VCO phase $\theta_2(0)$ and filter state $x(0)$. However, the phase error $\theta_\Delta$ may substantially increase during the acquisition process. In order to consider the property of the model to synchronize without undesired growth of the phase error $\theta_\Delta$, a lock-in range concept was introduced in \citep{Gardner-1966}:
``{\it{If, for some reason,  the frequency difference between input and VCO
is less than the loop bandwidth, the loop will lock up almost instantaneously
without slipping cycles. The maximum frequency difference for which
this fast acquisition is possible is called the lock-in frequency}}''.
The lock-in range concept is widely used in engineering literature on the PLL-based circuits study (see, e.g., \citep{Stensby-1997,KiharaOE-2002,Kroupa-2003,Gardner-2005-book,Best-2007}).
\begin{figure*}[!htbp]
\includegraphics[width=\linewidth]{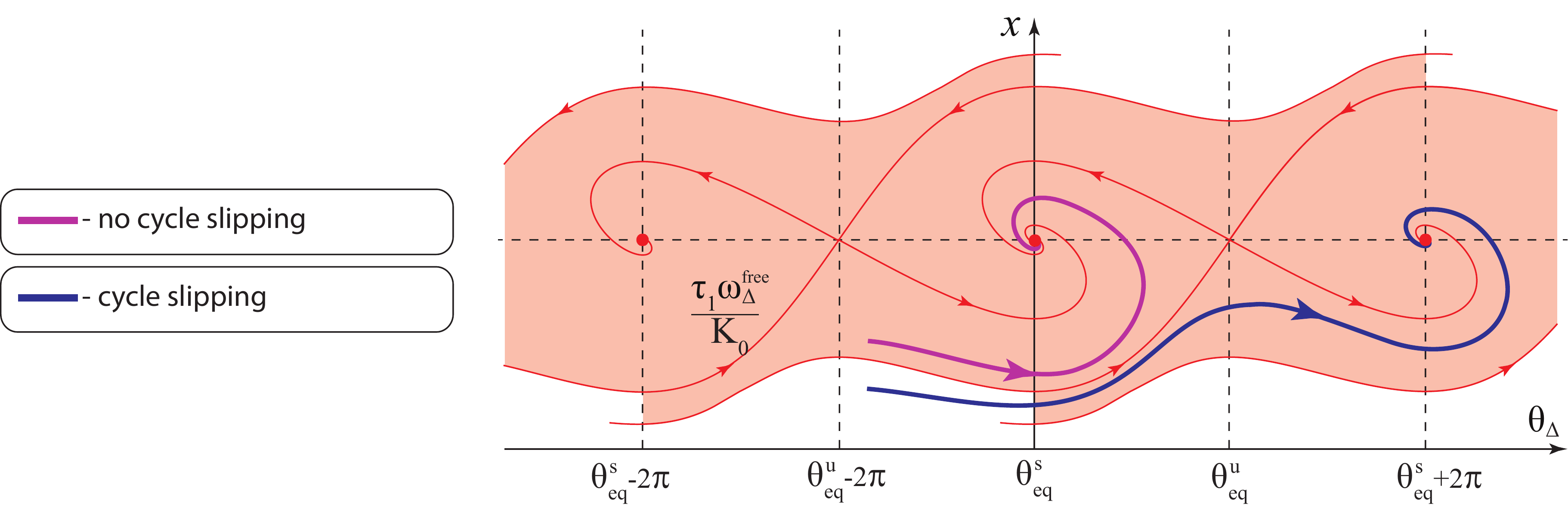}
\caption{The lock-in domain $D_{\rm lock-in}(\omega_\Delta^{\rm free})$ of (\ref{sys:PLLSys}).}
\label{ris:CycleSlipping}
\end{figure*}
Remark, that it is said that cycle slipping occurs if (see, e.g., \citep{AscheidM-1982,ErshovaL-1983,SmirnovaPU-2014}) $$\displaystyle \limsup_{t \rightarrow +\infty} \left|\theta_\Delta(0) - \theta_\Delta(t) \right| \geq 2\pi.$$
For (\ref{sys:PLLSys}) with fixed $\omega_\Delta^{\rm free}$ a domain of loop states for which the synchronization without cycle slipping occurs is called the lock-in domain $D_{\rm lock-in}(\omega_\Delta^{\rm free})$ (see Fig.~\ref{ris:CycleSlipping}).
However, in general, even for zero frequency deviation ($\omega_\Delta^{\text{free}}=0$)
and a sufficiently large initial state of filter ($x(0)$),
cycle slipping may take place, thus in 1979 Gardner wrote: \textit{``There is no natural way to define exactly any unique lock-in frequency''} and \textit{``despite its vague reality, lock-in range is a useful concept''} \citep{Gardner-1979-book}.
To overcome the stated problem, in \citep{KuznetsovLYY-2015-IFAC-Ranges,LeonovKYY-2015-TCAS} the rigorous mathematical definition of a lock-in range is suggested:
\begin{definition} \citep{KuznetsovLYY-2015-IFAC-Ranges,LeonovKYY-2015-TCAS}
\it The lock-in range of model (\ref{sys:PLLSys}) is a range $\left[0, \omega_l\right)$ such that for each frequency deviation $\left|\omega_\Delta^{\rm free}\right| \in \left[0, \omega_l\right)$ the model (\ref{sys:PLLSys}) is globally asymptotically stable and the following domain $$D_{\rm lock-in}\left((-\omega_l, \omega_l)\right) = \underset{\left|\omega_\Delta^{\rm free}\right|<\omega_l}{\bigcap} D_{\rm lock-in}(\omega_\Delta^{\rm free})$$ contains all corresponding equilibria $\left(\theta^s_{eq}, x_{eq}(\omega_\Delta^{\rm free})\right).$
\label{def:LockIn}
\end{definition}
For model (\ref{sys:PLLSys}) each lock-in domain from intersection $\underset{\left|\omega_\Delta^{\rm free}\right|<\omega_l}{\bigcap} D_{\rm lock-in}(\omega_\Delta^{\rm free})$ is bounded by the separatrices of saddle equilibria $\left(\theta^u_{eq}, x_{eq}(\omega_\Delta^{\rm free})\right)$ and vertical lines $\theta_\Delta = \theta^s_{eq} \pm 2\pi$.
Thus, the behavior of separatrices on the phase plane is the key to the lock-in range study (see Fig.~\ref{ris:LockinPhasePlane}).


\section{Phase plane analysis for the lock-in range estimation}
\label{sec:PhasePlaneAnalysis}
Consider an approach to the lock-in range computation of (\ref{sys:PLLSys}), based on the phase plane analysis.
To compute the lock-in range of (\ref{sys:PLLSys}) we need to consider the behavior of the lower separatrix $Q(\theta_\Delta, \omega_\Delta^{\rm free})$, which tends to the saddle point $\left(\theta_{eq}^u, x_{eq}(\omega_\Delta^{\rm free})\right) = \left(\pi, \frac{\omega_\Delta^{\rm free} \tau_1}{K_0}\right)$ as $t \rightarrow +\infty$ (by the symmetry of the lower and the upper half-planes, the consideration of the upper separatrix is also possible).

\begin{figure}[!htbp]
\includegraphics[width=0.9\linewidth]{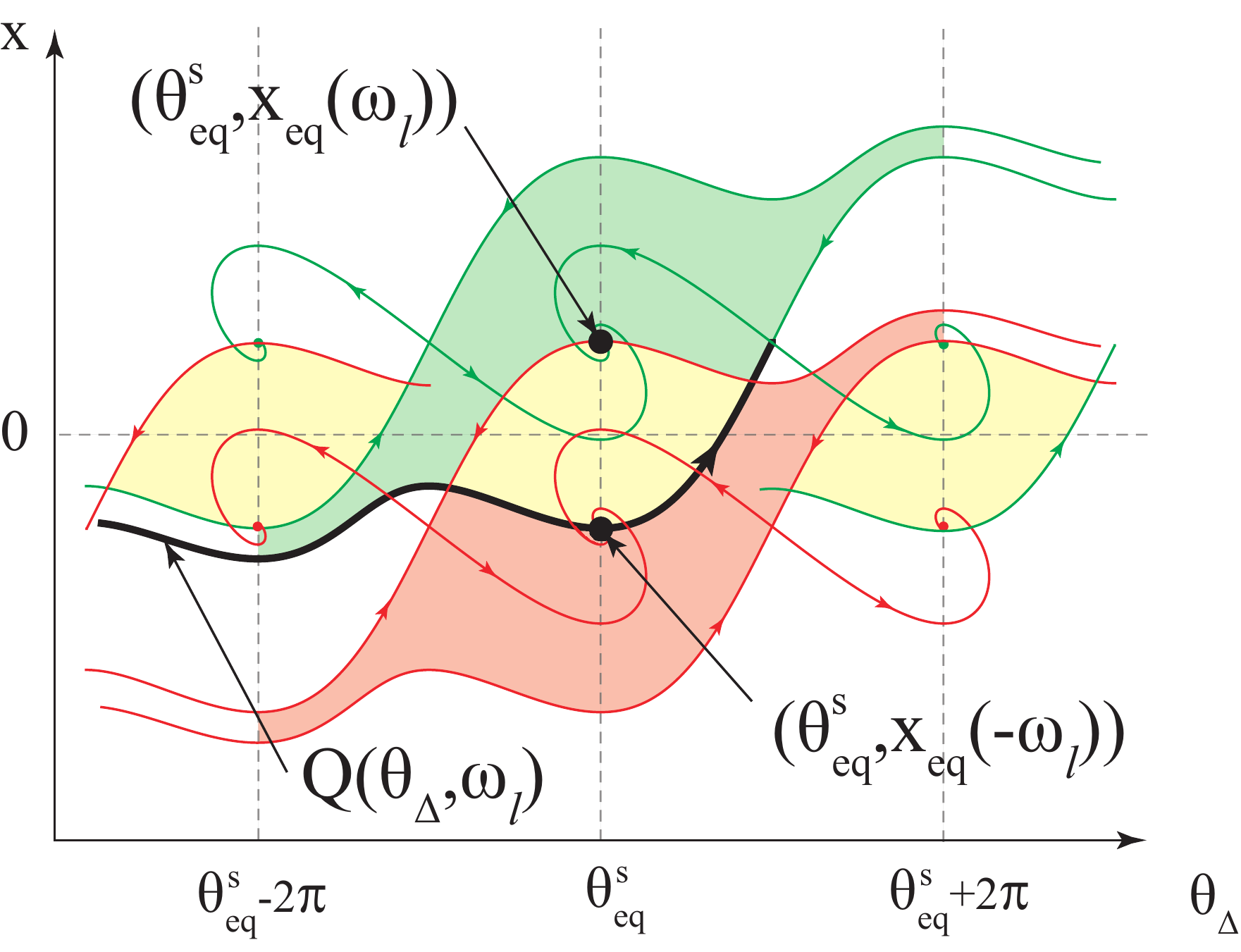}
\caption{The domain $D_{\rm lock-in}\left((-\omega_l, \omega_l)\right)$ of (\ref{sys:PLLSys}).}
\label{ris:LockinPhasePlane}
\end{figure}

The parameter $\omega_\Delta^{\rm free}$ shifts the phase plane vertically. To check this, we use a linear transformation $x \rightarrow x + \frac{\omega_\Delta^{\rm free} \tau_1}{K_0}$. Thus, to compute the lock-in range of (\ref{sys:PLLSys}), we need to find $\omega_\Delta^{\rm free} = \omega_l$ (where $\omega_l$ is called a lock-in frequency) such that (see Fig.~\ref{ris:LockinPhasePlane})
\begin{equation}
x_{eq}(-\omega_l) = Q(\theta^s_{eq}, \omega_l).
\label{eq:LockinRelation}
\end{equation}

By (\ref{eq:LockinRelation}), we obtain an exact formula for the lock-in frequency $\omega_l$:
\begin{align}
&-\frac{\omega_l}{K_0 / \tau_1} = \frac{\omega_l}{K_0 / \tau_1} + Q(\theta^s_{eq}, 0). \nonumber \\
&\omega_l = -\frac{K_0 Q(\theta^s_{eq}, 0)}{2 \tau_1},
\label{eq:LockinRelationNew}
\end{align}
\begin{figure}[!htbp]
\centering
\includegraphics[width=0.9\textwidth]{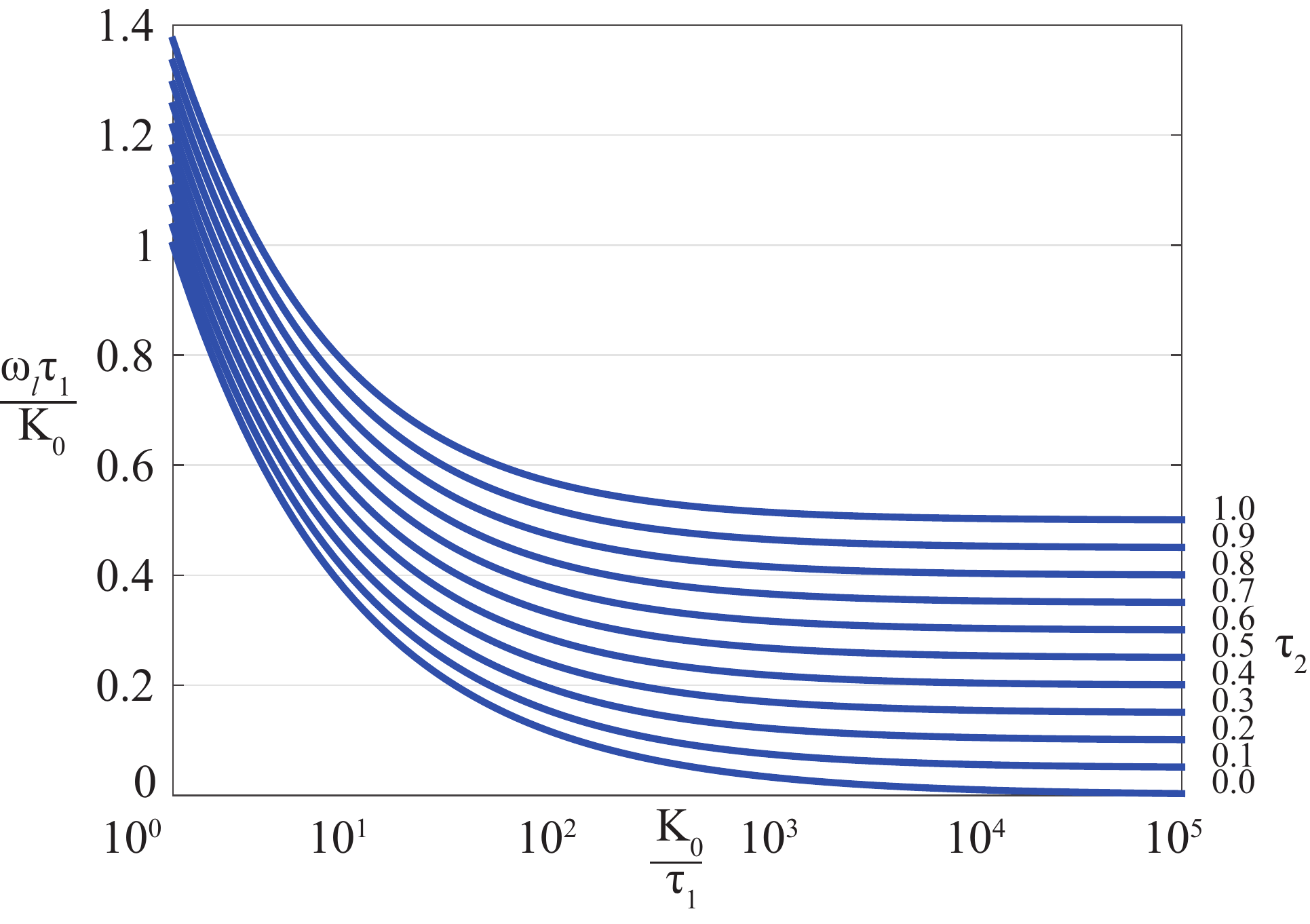}
\caption{Values of $\frac{\omega_l}{K_0 / \tau_1}$ for various $K_0$, $\tau_1$, $\tau_2$.}
\label{ris:Lock-inDiagram}
\end{figure}

Numerical simulations are used to compute the lock-in range of (\ref{sys:PLLSys}) applying (\ref{eq:LockinRelationNew}). The separatrix $Q(\theta_\Delta, 0)$ is numerically integrated and the corresponding $\omega_l$ is approximated. The obtained numerical results can be illustrated by a diagram (see Fig.~\ref{ris:Lock-inDiagram})\footnote{These results submitted to IFAC PSYCO 2016}.

Note that (\ref{sys:PLLSys}) depends on the value of two coefficients $\frac{K_0}{\tau_1}$ and $\tau_2$. In Fig.~\ref{ris:Lock-inDiagram}, choosing X-axis as $\frac{K_0}{\tau_1}$, we can plot a single curve for every fixed value of $\tau_2$. The results of numerical simulations show that for sufficiently large $\frac{K_0}{\tau_1}$, the value of $\omega_l$ grows almost proportionally to $\frac{K_0}{\tau_1}$. Hence, $\frac{\omega_l \tau_1}{K_0}$ is almost constant for sufficiently large $\frac{K_0}{\tau_1}$ and in Fig.~\ref{ris:Lock-inDiagram} the Y-axis can be chosen as $\frac{\omega_\Delta^{\rm free} \tau_1}{K_0}$.

\begin{figure}[!htbp]
\centering
\includegraphics[width=0.9\textwidth]{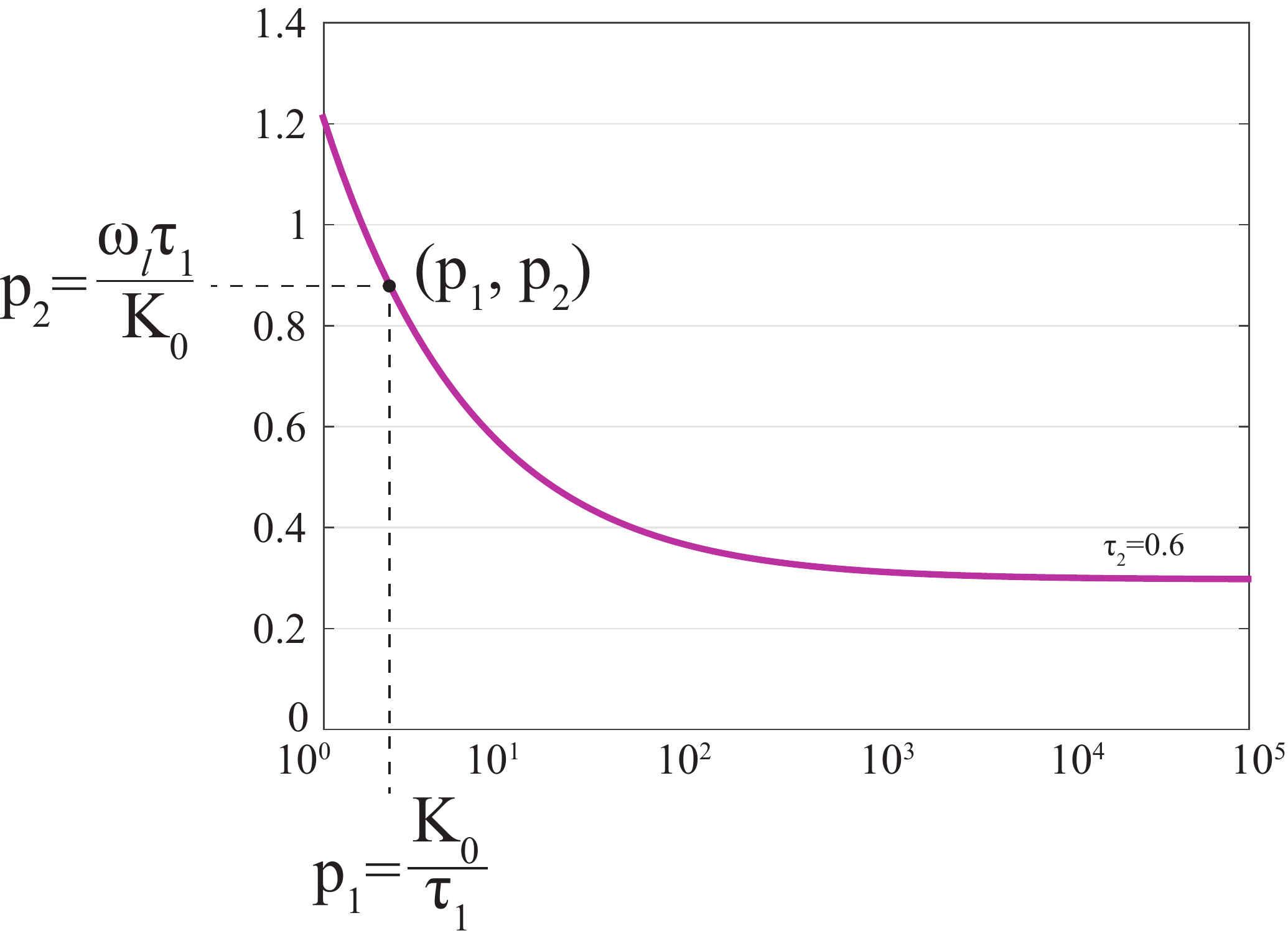}
\caption{The lock-in frequency calculation: $\omega_l=p_1p_2$.}
\label{ris:Lock-inDiagramExample}
\end{figure}

To obtain the lock-in frequency $\omega_l$ for fixed $\tau_1$, $\tau_2$, and $K_0$ using Fig.~\ref{ris:Lock-inDiagram}, we consider the curve corresponding to the chosen $\tau_2$. Next, for X-value equal $\frac{K_0}{\tau_1}$ we get the Y-value of the curve. Finally, we multiply the Y-value by $\frac{K_0}{\tau_1}$ (see Fig.~\ref{ris:Lock-inDiagramExample}).

Consider an analytical approach to the lock-in range estimation.
Main stages of the approach are presented in Subsection \ref{subsec:AnalyticalApproach}.
\subsection{Analytical approach to the lock-in range estimation}
\label{subsec:AnalyticalApproach}
Consider an active PI filter with small parameter $0 < \frac{\tau_2}{\tau_1} \ll 1$
(see, e.g., \citep{AlexandrovKLS-2014-ICUMT}). The consideration of (\ref{sys:PLLSys}) with such active PI filter allows us to estimate the lower separatrix $Q(\theta_\Delta, 0)$ and the lock-in range. For this purpose the approximations of separatrix $Q(\theta_\Delta, 0)$ in interval $0 \leq \theta_\Delta < \pi$ are used.

The separatrix $Q(\theta_\Delta, 0)$, which is a solution of (\ref{sys:PLLSys}), can be expanded in a Taylor series in variable $\tau_2/\tau_1$ (since the parameter $\tau_2/\tau_1$ is considered as a variable, the separatrix $Q(\theta_\Delta, 0) = Q(\theta_\Delta, 0, \tau_2/\tau_1)$ depends on it).
The first-order approximation of the lower separatrix $Q(\theta_\Delta, 0, \tau_2/\tau_1)$ has the form
\begin{align}
&\hat{Q}_{1}(\theta_\Delta, 0, \tau_2/\tau_1) = -2\sqrt{K_0 / \tau_1}\cos\frac{\theta_\Delta}{2} - \nonumber\\
&- \frac{\tau_2}{\tau_1}\frac{K_0 \left(\frac{2}{3} - \sin \frac{\theta_\Delta}{2} - \frac{1}{3} \sin \frac{3\theta_\Delta}{2}\right)}{\cos\frac{\theta_\Delta}{2}}.
\label{sys:SFirstApprox}
\end{align}
The second-order approximation of $Q(\theta_\Delta, 0, \tau_2/\tau_1)$ has the form
\begin{align}
&\hat{Q}_{2}(\theta_\Delta, 0, \tau_2/\tau_1) = -2\sqrt{K_0/ \tau_1}\cos\frac{\theta_\Delta}{2} - \frac{\tau_2}{\tau_1}\frac{K_0 \left(\frac{2}{3} - \sin \frac{\theta_\Delta}{2} - \frac{1}{3} \sin \frac{3\theta_\Delta}{2}\right)}{\cos\frac{\theta_\Delta}{2}} - \nonumber\\
& - \left(\frac{\tau_2}{\tau_1}\right)^2 \frac{K_0^2(6\frac{1}{2} - 4\ln 2)}{6\sqrt{K_0 /\tau_1}\cos\frac{\theta_\Delta}{2}} + \frac{K_0^2 \left(\frac{2}{3} - \sin \frac{\theta_\Delta}{2} - \frac{1}{3} \sin \frac{3\theta_\Delta}{2}\right)^2}{4\sqrt{K_0 /\tau_1} \cos^3\frac{\theta_\Delta}{2}}+ \nonumber \\
& + \left(\frac{\tau_2}{\tau_1}\right)^2\frac{K_0^2\left(8 \sin(\frac{\theta_\Delta}{2}) - 4\ln \Big|\sin \frac{\theta_\Delta}{2} + 1\Big|\right)}{6\sqrt{K_0 /\tau_1}\cos\frac{\theta_\Delta}{2}} + \left(\frac{\tau_2}{\tau_1}\right)^2\frac{K_0^2\left( \frac{1}{2}\cos 2\theta_\Delta + 2\cos \theta_\Delta\right)}{6\sqrt{K_0 /\tau_1}\cos\frac{\theta_\Delta}{2}}.
\label{sys:SSecondApprox}
\end{align}

For approximations (\ref{sys:SFirstApprox}), (\ref{sys:SSecondApprox}) of separatrix $Q(\theta_\Delta, 0, \tau_2/\tau_1)$ the following relations are valid:
\begin{align*}
&Q(\theta_\Delta, 0, \tau_2/\tau_1) = \hat{Q}_{1}(\theta_\Delta, 0, \tau_2/\tau_1) + O\left(\left(\tau_2/\tau_1\right)^2\right), \\
&Q(\theta_\Delta, 0, \tau_2/\tau_1) = \hat{Q}_{2}(\theta_\Delta, 0, \tau_2/\tau_1) + O\left(\left(\tau_2/\tau_1\right)^3\right).
\end{align*}
For $\theta_\Delta = \theta^s_{eq}$ the relations (\ref{sys:SFirstApprox}), (\ref{sys:SSecondApprox}) take the following values:
\begin{align*}
&\hat{Q}_1(\theta^s_{eq}, 0, \tau_2/\tau_1) = - 2 \sqrt{K_0 / \tau_1} - \frac{2K_0 \tau_2}{3 \tau_1}, \\
&\hat{Q}_2(\theta^s_{eq}, 0, \tau_2/\tau_1) = - 2 \sqrt{K_0 / \tau_1} - \frac{2K_0 \tau_2}{3 \tau_1} - \\
& -\frac{K_0 \tau_2^2 (5 - 6\ln 2)}{9 \tau_1}\sqrt{K_0 / \tau_1}.
\end{align*}
Using relation (\ref{eq:LockinRelationNew}) the lock-in frequency $\omega_l$ is approximated as follows:
\begin{align}
&\omega_l = \frac{K_0 \sqrt{K_0 / \tau_1}}{\tau_1} + \frac{K_0^2\tau_2}{3 \tau_1^2} + O\left(\left(\tau_2 / \tau_1\right)^2\right),
\label{eq:LockinFirstApprox}
\end{align}
\begin{align}
&\omega_l = \frac{K_0 \sqrt{K_0 / \tau_1}}{\tau_1} + \frac{K_0^2\tau_2}{3 \tau_1^2} + \nonumber \\
& + \frac{K_0^2 \tau_2^2(5 - 6\ln 2)}{18 \tau_1^2}\sqrt{K_0 / \tau_1} + O\left(\left(\tau_2 / \tau_1\right)^3\right).
\label{eq:LockinSecondApprox}
\end{align}
\begin{figure}[!htbp]
\centering
\includegraphics[width=0.7\textwidth]{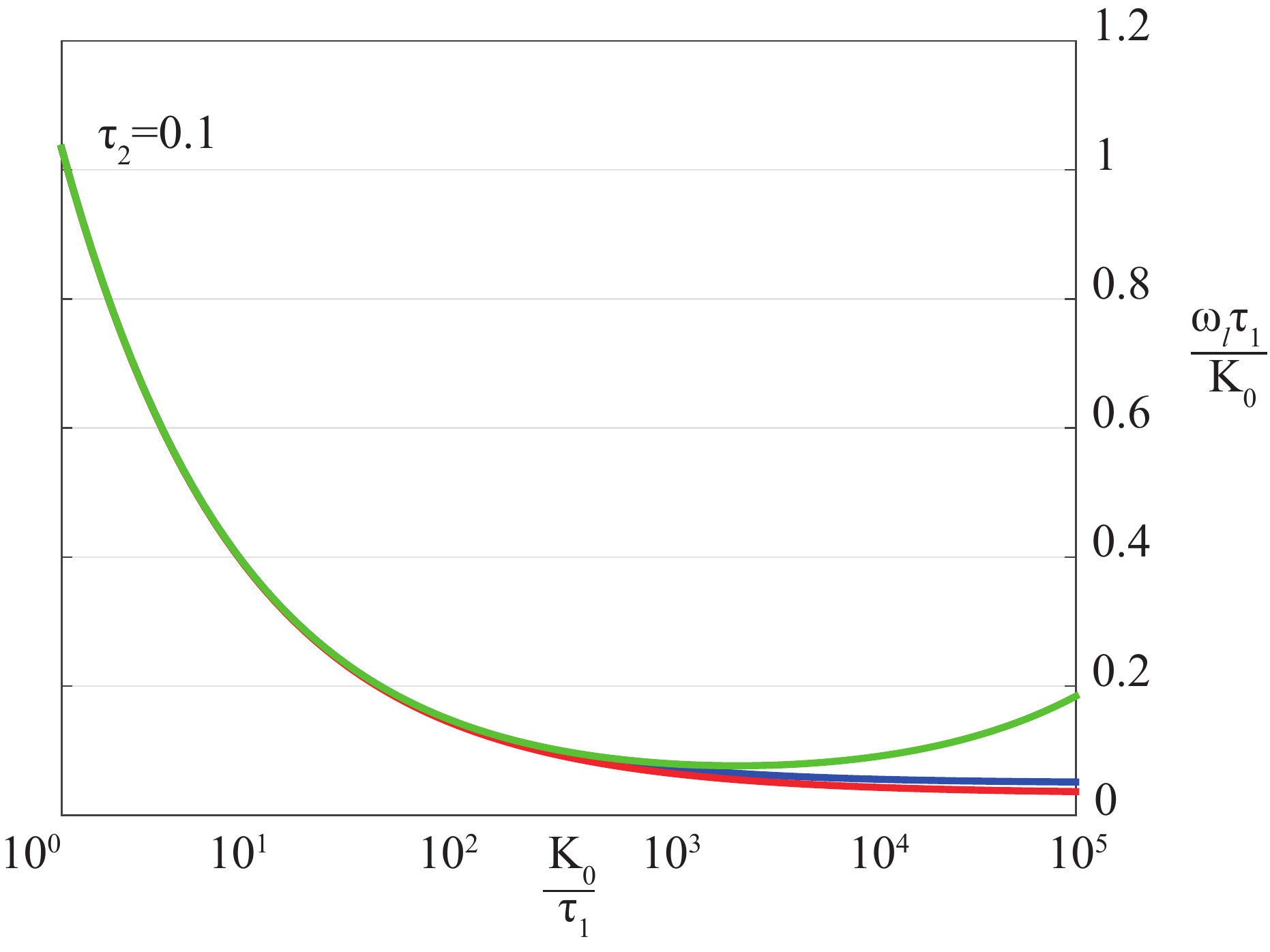}
\caption{Estimates on $\frac{\omega_l}{K_0 / \tau_1}$ for various $K_0$, $\tau_1$.}
\label{ris:ApproxExample}
\end{figure}

For fixed $\tau_2 = 0.1$ the three curves are shown in Fig.~\ref{ris:ApproxExample}. The values of $\omega_l$ (the blue curve, which is obtained numerically using relation (\ref{eq:LockinRelationNew})) are estimated from below by (\ref{eq:LockinFirstApprox}) and from above by (\ref{eq:LockinSecondApprox}) (the red and green curves correspondingly).
Since the lock-in frequency $\omega_l$ is approximated under the condition of small parameter $\tau_2/\tau_1$, the estimates (\ref{eq:LockinFirstApprox}) and (\ref{eq:LockinSecondApprox}) give less precise result in the case of large $K_0/\tau_1$.

\subsection{The pull-out frequency and lock-in range}
\label{subsec:PullOutLockIn}
An another characteristic related to the cycle slipping effect is the pull-out frequency $\omega_{\rm po}$ (see, e.g., \citep{Gardner-1979-book, Stensby-1997, Kroupa-2003}. In \citep{Gardner-2005-book} the pull-out frequency is defined as a frequency-step limit, \textit{``below which the loop does not skip cycles but remains in lock''}. However, in general case of Filter (see, e.g., \citep{PinheiroP-2014,BanerjeeS-2008}) the pull-out frequency may depend on the value of $\omega_\Delta^{\rm free}$.

\begin{figure}[!htbp]
\includegraphics[width=0.8\linewidth]{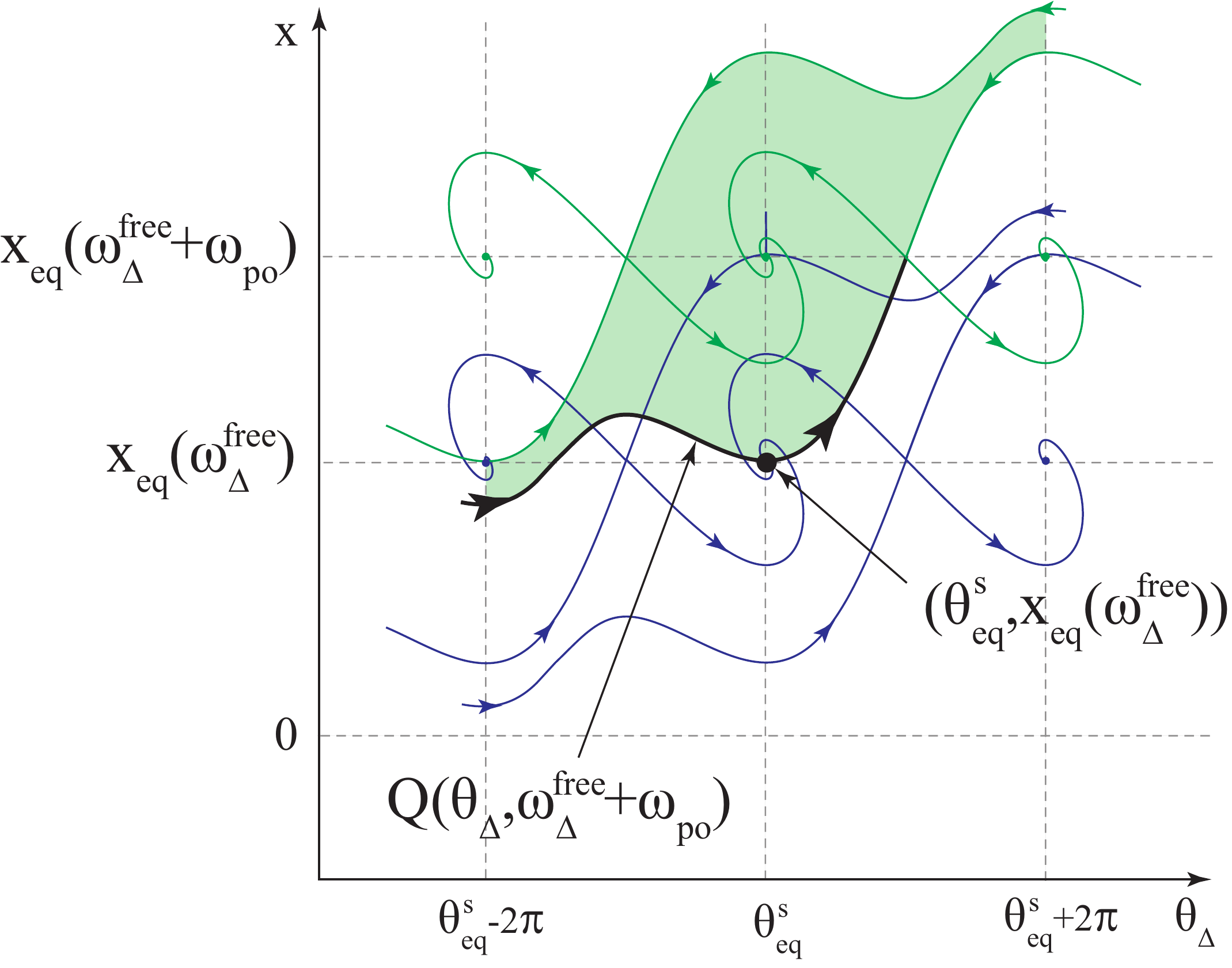}
\caption{The frequency step of (\ref{sys:PLLSys}) equals to pull-out frequency $\omega_{po}$.}
\label{ris:PullOutPhasePlane}
\end{figure}

However, in the case of active PI filter, the pull-out frequency can be defined and approximated (see, e.g., \citep{Gardner-1979-book,HuqueS-2013}), since the parameter $\omega_\Delta^{\rm free}$ only shifts the phase plane vertically. The pull-out frequency can be found as follows (see Fig.~\ref{ris:PullOutPhasePlane}):
\begin{align}
&x_{eq}(\omega_\Delta^{\rm free}) = Q(\theta^s_{eq}, \omega_\Delta^{\rm free} + \omega_{po}), \nonumber \\
&\frac{\omega_\Delta^{\rm free}}{K_0 / \tau_1} = \frac{\omega_\Delta^{\rm free} + \omega_{po}}{K_0 / \tau_1} + Q(\theta^s_{eq}, 0). \nonumber \\
&\omega_{po} = -\frac{K_0 Q(\theta^s_{eq}, 0)}{\tau_1} = 2 \omega_l.
\label{eq:PullOutRelation}
\end{align}

\begin{figure*}[!htbp]
\centering
\includegraphics[width=0.9\textwidth]{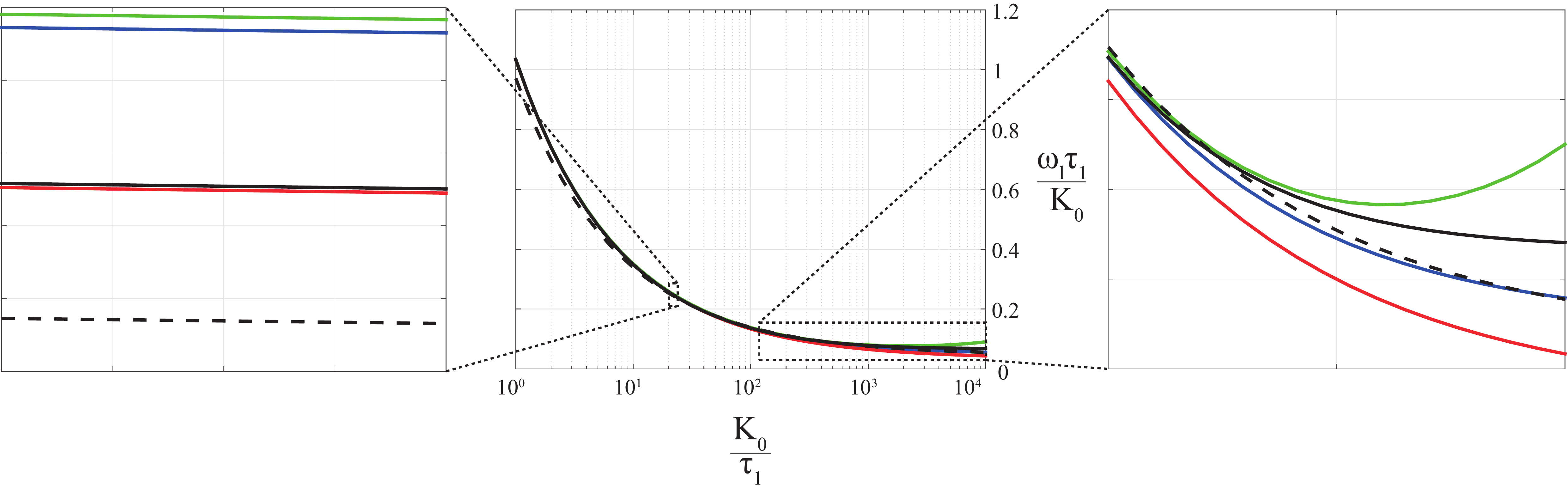}
\caption{Comparison of the pull-out frequency estimates.}
\label{ris:Comparison}
\end{figure*}

In Fig.~\ref{ris:Comparison} the estimates from \citep{Gardner-1979-book,HuqueS-2013} are compared with estimates based on (\ref{eq:LockinFirstApprox}) and (\ref{eq:LockinSecondApprox}).
The pull-out frequency estimate, which is obtained according to Fig.~\ref{ris:Lock-inDiagram} and (\ref{eq:PullOutRelation}), is drawn in blue color.
Analytical estimates based on (\ref{eq:LockinFirstApprox}), (\ref{eq:LockinSecondApprox}), and (\ref{eq:PullOutRelation}) are drawn in red and green colors correspondingly.
The black curve is the estimate of the pull-out frequency from \citep{HuqueS-2013}. The dashed curve corresponds to the empirical estimate
\begin{equation}
\omega_{po} \approx 1.85\left(\frac{1}{2} + \frac{\tau_1}{K_0 \tau_2^2}\right),
\label{rel:ViterbiEmp}
\end{equation}
presented in \citep{Gardner-1979-book}.

For $K_0/\tau_1$ not very large the relation (\ref{eq:LockinSecondApprox}) is the most precise estimate compared to the presented ones.

\section{Conclusion}
In the present work models of the PLL-based circuits in the signal's phase space are described. The lock-in range of PLL-based circuits with sinusoidal PD characteristic and active PI filter is considered. The rigorous definition of the lock-in range is discussed, and relation (\ref{eq:LockinRelationNew}) for the lock-in range computation is derived. For the lock-in range estimation two approaches -- numerical and analytical -- are presented. The methods are based on the integration of phase trajectories. In Subsection \ref{subsec:AnalyticalApproach} the numerical estimates are verified by analytical estimates, which are obtained under the condition of small parameter.

\appendix

\section{The lock-in range estimation for small parameter of the loop filter.}
\label{sec:SmallParProof}

Let us write out (\ref{sys:PLLSys}) in a different form with $a = \frac{\tau_2}{\tau_1}$ and $b = \frac{1}{\tau_1}$:
\begin{equation}
\begin{cases}
\dot{x} = \sin(\theta_\Delta), \\
\dot{\theta}_\Delta = \omega_\Delta^{\rm free} - b K_0 x - a K_0 \sin(\theta_\Delta).
\end{cases}
\label{app:PLLSysOrig}
\end{equation}

Consider the following system, which is equivalent to (\ref{app:PLLSysOrig}):
\begin{equation}
\begin{cases}
\dot{\theta}_\Delta = y,\\
\dot{y} = - aK_0\cos(\theta_\Delta)y - {b}K_0\sin(\theta_\Delta),
\label{eq:AppPLLSys}
\end{cases}
\end{equation}
where $y = \omega_\Delta^{\rm free} - b K_0 x - a K_0 \sin(\theta_\Delta)$.

In virtue of $2\pi$-periodicity of (\ref{eq:AppPLLSys}) in variable $\theta_\Delta$, phase trajectories of (\ref{eq:AppPLLSys}) coincides for each interval $\theta_\Delta \in \left(-\pi + 2\pi k, \pi + 2\pi k\right]$, $k \in \mathbb{Z}$. Thus, one can study (\ref{eq:AppPLLSys}) in interval $\theta_\Delta \in \left(-\pi, \pi\right]$ only.

Let us find equilibria of (\ref{eq:AppPLLSys}) from the following system of equations:
\begin{equation}
\begin{cases}
\sin(\theta_{eq}) = 0, \\
a K_0 \cos(\theta_{eq}) y_{eq} = 0.
\end{cases}
\nonumber
\end{equation}

In interval $\theta_\Delta \in \left(-\pi, \pi\right]$ there exist two equilibria $\left(\theta^s_{eq}, y_{eq}\right) = (0; 0)$ and $\left(\theta^u_{eq}, y_{eq}\right) = (\pi; 0)$. To define type of the equilibria points let us write out corresponding characteristic polynomials and find the eigenvalues:
\begin{equation}
\begin{array}{ll}
\text{equilibrium }(0; 0): & \lambda^2 + aK_0\lambda + b K_0 = 0; \\
\\
\lambda_{1,2} = \frac{-aK_0 \pm \sqrt{(aK_0)^2 - 4b K_0}}{2}, & (aK_0)^2 - 4b K_0>0; \\
\lambda_{1} = \lambda_{2} = \frac{-aK_0}{2}, \hskip0.2cm & (aK_0)^2 - 4b K_0=0;\\
\lambda_{1,2} = \frac{-aK_0 \pm i \sqrt{4b K_0 - (aK_0)^2}}{2}, & (aK_0)^2 - 4b K_0<0;\\
\end{array} \nonumber
\end{equation}
\\
\begin{equation}
\begin{array}{ll}
\text{equilibrium }(\pi; 0): & \lambda^2 - aK_0\lambda - b K_0 = 0;\\
\\
\lambda_{1,2} = \frac{aK_0 \pm \sqrt{(aK_0)^2 + 4b K_0}}{2}. &
\end{array} \nonumber
\end{equation}
Thus, equilibrium $\left(\theta^s_{eq}, y_{eq}\right)$ is a stable node, a stable degenerated node, or a stable focus (that depends on the sign of $(aK_0)^2 - 4b K_0$). Equilibrium $\left(\theta^u_{eq}, y_{eq}\right)$ is a saddle point for all $a > 0$, $b > 0$, $K_0 > 0$.
Moreover, in virtue of periodicity each equilibrium $\left(\theta^u_{eq} + 2\pi k, y_{eq}\right)$ is a saddle point, and each equilibrium $\left(\theta^s_{eq} + 2\pi k, y_{eq}\right)$ is a stable equilibrium of the same type as $\left(\theta^s_{eq}, y_{eq}\right)$.
Note also that equilibria $\left(\theta_{eq}, y_{eq}\right)$ of (\ref{eq:AppPLLSys}) and corresponding equilibria $\left(\theta_{eq}, x_{eq}\right)$ of (\ref{app:PLLSysOrig}) are of the same type, and related as follows:
\begin{equation}
\left(\theta_{eq}, y_{eq}\right) = \left(\theta_{eq}, \omega_\Delta^{\rm free} - b K_0 x_{eq}\right).
\nonumber
\end{equation}

Let us consider the following differential equation:
\begin{equation}
y'(\theta_\Delta) = \displaystyle -\frac{a K_0\cos(\theta_\Delta)y(\theta_\Delta) + {b}K_0\sin(\theta_\Delta)}{y(\theta_\Delta)}.
\label{eq:PLLEqPhasePlane}
\end{equation}
The right side of equation (\ref{eq:PLLEqPhasePlane}) is discontinuous in each point of line $y=0$. This line is an isocline line of vertical angular inclination of (\ref{eq:PLLEqPhasePlane}) \citep{BarbashinT-1969}. Equation (\ref{eq:PLLEqPhasePlane}) is equivalent to (\ref{eq:AppPLLSys}) in the upper and the lower open half planes of the phase plane.

Let the solutions $y(\theta_\Delta, a)$ of equation (\ref{eq:PLLEqPhasePlane}) be considered as functions of two variables $\theta_\Delta$, $a$.
Consider the solution of differential equation (\ref{eq:PLLEqPhasePlane}), which range of values lies in the upper open half plane of its phase plane. Right side of equation (\ref{eq:PLLEqPhasePlane}) in the upper open half plane is function of class $C^m$ for $m$ arbitrary large. Solutions of the Cauchy problem with initial conditions $x = x_0$, $y = y_0$ (which solutions are on the upper half plane) are also of class $C^m$ on their domain of existence for $m$ arbitrary large \citep{Hartman-1964}.

Let us study the separatrix $S(\theta_\Delta, a)$ in interval $\left[0, \pi\right)$, which tends to saddle point $(\theta^u_{eq}; x_{eq}) = (\pi; 0)$ and is situated in its second quadrant. Separatrix $S(\theta_\Delta, a)$ is the solution of the corresponding Cauchy problem for equation (\ref{eq:PLLEqPhasePlane}). The separatrix $S(\theta_\Delta, a)$ is of class $C^m$ on its domain of existence for $m$ arbitrary large.

Consider separatrix $S(\theta_\Delta, a)$ as a Taylor series in variable $a$ in the neighborhood of $a_0 = 0$:
\begin{equation}
S(\theta_\Delta, a) = S(\theta_\Delta, 0) + a\frac{\partial S(\theta_\Delta, a)}{\partial a}\Big|_{a=0} + \dots + \frac{a^n}{n!}\frac{\partial^n S(\theta_\Delta, a)}{\partial a^n}\Big|_{a=0} + \dots.
\label{eq:TailorS}
\end{equation}

Let us denote
\begin{align*}
&S_0(\theta_\Delta) = S(\theta_\Delta, 0), \\
&S_{i}(\theta_\Delta) = \frac{1}{i!}\frac{\partial^{i} S(\theta_\Delta, a)}{\partial a^{i}}\Big|_{a=0}, \hskip0.5cm i \geq 1.
\end{align*}
$\hat{S}_n(\theta_\Delta, a)$ as the $n$-th approximation of $S(\theta_\Delta, a)$ in variable $a$:
\begin{equation}
\hat{S}_n(\theta_\Delta, a) = S(\theta_\Delta, 0) + \sum \limits_{i=1}^{n} a^i S_{i}(\theta_\Delta). \nonumber
\end{equation}
The Taylor remainder is denoted as follows:
\begin{equation}
\tilde{S}_{n}(\theta_\Delta, a) = \sum \limits_{i=n+1}^{+\infty} a^i S_{i}(\theta_\Delta).
\end{equation}
For the convergent Taylor series its remainder $\tilde{S}_{n}(\theta_\Delta, a) = O(a^{n+1})$ for each point $x_0$ of interval $\left[0, \pi\right)$.

Separatrix $S(\theta_\Delta, a)$ satisfies the following relation, which follows from (\ref{eq:PLLEqPhasePlane}):
\begin{equation}
S(\theta_\Delta, a)S^{\prime}(\theta_\Delta, a) = -a K_0\cos(\theta_\Delta)S(\theta_\Delta, a) - {b}K_0\sin(\theta_\Delta).
\nonumber
\end{equation}
\begin{equation}
\int \limits_{\theta_\Delta}^{\pi}S(s, a)dS(s, a) = -a K_0\int \limits_{\theta_\Delta}^{\pi}\cos(s)S(s, a)ds - {b}K_0\int \limits_{\theta_\Delta}^{\pi}\sin(s)ds.
\nonumber
\end{equation}
\begin{equation}
\frac{1}{2}\displaystyle \lim_{s \rightarrow \pi-0} S^2(s,a) - \frac{1}{2}S^2(\theta_\Delta,a)= -a K_0\int \limits_{\theta_\Delta}^{\pi}\cos(s)S(s, a)ds - {b}K_0\int \limits_{\theta_\Delta}^{\pi}\sin(s)ds.
\label{eq:SInEqInt}
\end{equation}
Let us represent $S(\theta_\Delta, a)$ as Taylor series (\ref{eq:TailorS}) in relation (\ref{eq:SInEqInt}).
\begin{align*}
&-\frac{1}{2}\left(S_0(\theta_\Delta) + a S_1(\theta_\Delta) + a^2 S_2(\theta_\Delta) + \tilde{S}_2(\theta_\Delta, a)\right)^2 =  - {b}K_0\int \limits_{\theta_\Delta}^{\pi}\sin(s)ds - \\
&-a K_0\int \limits_{\theta_\Delta}^{\pi}\cos(s)\left(S_0(s) + S_1(s) + S_2(s) + \tilde{S}_2(s, a)\right)ds.
\end{align*}
\begin{align*}
&-\frac{1}{2}S_0^2(\theta_\Delta) - a S_0(\theta_\Delta) S_1(\theta_\Delta) - \frac{1}{2}a^2 S_1^2(\theta_\Delta) - a^2 S_0(\theta_\Delta) S_2(\theta_\Delta) + O(a^3) = \\
& - {b}K_0\int \limits_{\theta_\Delta}^{\pi}\sin(s)ds - a K_0\int \limits_{\theta_\Delta}^{\pi}\cos(s)S_0(s) - a^2 K_0\int \limits_{\theta_\Delta}^{\pi}\cos(s)S_1(s)ds - O(a^3).
\end{align*}
\begin{align}
\nonumber
&-\frac{1}{2}S_0^2(\theta_\Delta) - a S_0(\theta_\Delta) S_1(\theta_\Delta) - a^2 \left( \frac{1}{2} S_1^2(\theta_\Delta) + S_0(\theta_\Delta) S_2(\theta_\Delta)\right) + O(a^3) = \\
& - {b}K_0\int \limits_{\theta_\Delta}^{\pi}\sin(s)ds - a K_0\int \limits_{\theta_\Delta}^{\pi}\cos(s)S_0(s) - a^2 K_0\int \limits_{\theta_\Delta}^{\pi}\cos(s)S_1(s)ds + O(a^3).
\label{eq:TaylorInEq}
\end{align}
Let us write out the corresponding members of (\ref{eq:TaylorInEq}) for each $a^n$, $n = 0, 1, 2$.\\
For $a^0$:
\begin{equation}
\frac{1}{2}S_0^2(\theta_\Delta) = b K_0 \int \limits_{\theta_\Delta}^{\pi} \sin(s) ds. \label{eq:approxEps0}
\end{equation}
For $a^1$:
\begin{equation}
S_{0}(\theta_\Delta)S_{1}(\theta_\Delta) = K_0\int \limits_{\theta_\Delta}^{\pi}\cos(s)S_0(s)ds.\label{eq:approxEps1}
\end{equation}
For $a^2$:
\begin{equation}
S_{0}(\theta_\Delta)S_{2}(\theta_\Delta) + \frac{1}{2} S_1^2(\theta_\Delta) = K_0\int \limits_{\theta_\Delta}^{\pi}\cos(s)S_1(s)ds. \label{eq:approxEps2}
\end{equation}

Let us consequently find $S_{0}(\theta_\Delta)$, $S_{1}(\theta_\Delta)$, $S_{2}(\theta_\Delta)$ using relations (\ref{eq:approxEps0}), (\ref{eq:approxEps1}) and (\ref{eq:approxEps2}). Begin with evaluation of $S_{0}(\theta_\Delta)$:
\begin{align}
&\frac{1}{2}S_0^2(\theta_\Delta) = b K_0 \int \limits_{\theta_\Delta}^{\pi} \sin(s) ds = - b K_0 \cos(s) \Big|_x^\pi ds = \nonumber \\
& = b K_0 (1+\cos(\theta_\Delta)). \nonumber \\
&S_0(\theta_\Delta) = \sqrt{2 b K_0(1+\cos(\theta_\Delta))}. \label{rel:approxEps0}
\end{align}
According to (\ref{rel:approxEps0})
\begin{equation}
S_0(0) = 2\sqrt{b K_0}.
\label{rel0:approxEps0}
\end{equation}

Using equation (\ref{eq:approxEps1}) and relations (\ref{rel:approxEps0}) evaluate $S_{1}(\theta_\Delta)$:
\begin{equation}
S_{1}(\theta_\Delta) = \frac{K_0\int \limits_{\theta_\Delta}^{\pi}\cos(s)S_0(s)ds}{S_0(\theta_\Delta)}. \nonumber
\end{equation}
\begin{equation}
S_{1}(\theta_\Delta) =\displaystyle \frac{K_0\sqrt{2{b}K_0}\int
\limits_{\theta_\Delta}^{\pi} \cos(s)\sqrt{(1+\cos(s))}ds}{\sqrt{2{b}K_0(1+\cos(\theta_\Delta))}}. \nonumber
\end{equation}
Let us evaluate the integral  $$\int \limits_{\theta_\Delta}^{\pi}
\cos(s)\sqrt{(1+\cos(s))}ds$$ in the interval $\theta_\Delta \in \left[0; \pi\right)$ using the following substitutions:
\begin{equation}
\begin{array}{l}  \nonumber
u = 1+\cos(s), du = -\sin(s)ds\\
v = \sqrt{2-u}, dv = \displaystyle -\frac{du}{2\sqrt{2-u}}.
\end{array}
\end{equation}

\begin{align*}
&\int \limits_{\theta_\Delta}^{\pi} \cos(s)\sqrt{(1+\cos(s))}ds = \int \limits_{\theta_\Delta}^{\pi} \frac{\cos(s)\sin(s)}{\sqrt{(1-\cos(s))}}ds =\\
& = -\int\limits^{0}_{1+\cos(\theta_\Delta)} \displaystyle \frac{u-1}{\sqrt{2-u}}du = 2\int \limits^{\sqrt{2}}_{\sqrt{1-\cos(\theta_\Delta)}} (1-v^2)dv = \\
&= 2\left(\sqrt{2}-\sqrt{1-\cos(\theta_\Delta)}\right) - \frac{2}{3}\left(2\sqrt{2} - \sqrt{1-\cos(\theta_\Delta)}^3\right) = \\
&= -\frac{2}{3}(2+\cos(\theta_\Delta))\sqrt{1-\cos(\theta_\Delta)} + \frac{2\sqrt{2}}{3}.
\end{align*}
Hence, an expression for $S_{1}(\theta_\Delta)$ in interval $\theta_\Delta \in \left[0; \pi\right)$ is obtained:
\begin{equation}
S_{1}(\theta_\Delta) = \displaystyle \frac{K_0\sqrt{2{b}K_0}\left(\frac{2\sqrt{2}}{3} - \frac{2}{3}(2+\cos(\theta_\Delta))\sqrt{1-\cos(\theta_\Delta)}\right)}{\sqrt{2{b}K_0(1+\cos(\theta_\Delta))}}. \label{rel:approxEps1}
\end{equation}
Moreover,
\begin{equation}
S_1(0) = \frac{2K_0}{3}. \label{rel0:approxEps1}
\end{equation}
To shorten the further evaluation of $S_2(\theta_\Delta)$, write out $S_1(\theta_\Delta)$ in equivalent form (in interval $\theta_\Delta \in \left[0; \pi\right)$).
\begin{align*}
&S_{1}(\theta_\Delta) = \displaystyle \frac{K_0\sqrt{2{b}K_0}\left(\frac{2\sqrt{2}}{3} - \frac{2}{3}(2+\cos(\theta_\Delta))\sqrt{1-\cos(\theta_\Delta)}\right)}{\sqrt{2{b}K_0(1+\cos(\theta_\Delta))}} =
\end{align*}
\begin{align*}
&= \displaystyle \frac{K_0\sqrt{2{b}K_0}\left(\frac{2\sqrt{2}}{3} - \frac{2}{3}(2+\cos \theta_\Delta)\sqrt{2}\sin\frac{\theta_\Delta}{2}\right)}{\sqrt{2{b}K_0}\sqrt{2}\cos\frac{\theta_\Delta}{2}} = \displaystyle \frac{K_0\left(\frac{2}{3} - \frac{2}{3}(2+\cos \theta_\Delta)\sin\frac{\theta_\Delta}{2}\right)}{\cos\frac{\theta_\Delta}{2}} =
\end{align*}
\begin{align*}
&= \displaystyle \frac{K_0\left(\frac{2}{3} - \frac{2}{3}(3-2\sin^2\frac{\theta_\Delta}{2})\sin\frac{\theta_\Delta}{2}\right)}{\cos\frac{\theta_\Delta}{2}} = \displaystyle \frac{K_0\left(\frac{2}{3} - 2\sin\frac{\theta_\Delta}{2} - \frac{4}{3}\sin^3\frac{\theta_\Delta}{2}\right)}{\cos\frac{\theta_\Delta}{2}} =
\end{align*}
\begin{align*}
&= \frac{K_0 \left(\frac{2}{3} - \sin \frac{\theta_\Delta}{2} - \frac{1}{3} \sin \frac{3\theta_\Delta}{2}\right)}{\cos\frac{\theta_\Delta}{2}}.
\end{align*}

Let us evaluate $S_2(\theta_\Delta)$ using (\ref{eq:approxEps2}), (\ref{rel:approxEps0}) and (\ref{rel:approxEps1}).
\begin{align*}
&S_{2}(\theta_\Delta) = \frac{K_0\int \limits_{\theta_\Delta}^{\pi}\cos(s)S_1(s)ds - \frac{1}{2} S_1^2(\theta_\Delta)}{S_0(\theta_\Delta)} = \\
& = \frac{K_0^2\int \limits_{\theta_\Delta}^{\pi}\frac{\cos(s)\left(\frac{2}{3} - \sin \frac{s}{2} - \frac{1}{3} \sin \frac{3s}{2}\right)}{\cos\frac{s}{2}}ds}{2\sqrt{{b} K_0}\cos\frac{\theta_\Delta}{2}} - \frac{K_0^2 \left(\frac{2}{3} - \sin \frac{\theta_\Delta}{2} - \frac{1}{3} \sin \frac{3\theta_\Delta}{2}\right)^2}{4\sqrt{{b} K_0} \cos^3 \frac{\theta_\Delta}{2}}.
\end{align*}
Evaluate the integral $\displaystyle \int \limits_{\theta_\Delta}^{\pi}\frac{\cos(s)\left(\frac{2}{3} - \sin \frac{s}{2} - \frac{1}{3} \sin \frac{3s}{2}\right)}{\cos\frac{s}{2}}ds$:

I
\begin{align*}
&\int \limits_{\theta_\Delta}^{\pi}\frac{2}{3} \frac{\cos(s)}{\cos\frac{s}{2}} ds= \frac{2}{3}\int \limits_{\theta_\Delta}^{\pi} \frac{2\cos^2(\frac{s}{2}) - 1}{\cos\frac{s}{2}} ds= \frac{2}{3}\int \limits_{\theta_\Delta}^{\pi} \left(2\cos\frac{s}{2} - \frac{1}{\cos\frac{s}{2}}\right) ds= \\
&= \frac{8}{3} \left(\sin(\frac{s}{2})\right)\Big|_x^\pi - \frac{2}{3}\int \limits_{\theta_\Delta}^{\pi} \frac{1}{\cos\frac{s}{2}} ds =
\end{align*}

\begin{align*}
& u = \frac{s}{2}; \hskip0.5cm du = \frac{1}{2}ds
\end{align*}

\begin{align*}
&= \frac{8}{3} \left(\sin(\frac{s}{2})\right)\Big|_x^\pi - \frac{4}{3}\int \limits_{\frac{\theta_\Delta}{2}}^{\frac{\pi}{2}} \frac{1}{\cos u} du = \left(\frac{8}{3} \sin(\frac{s}{2}) - \frac{4}{3}\ln \Big|\operatorname{tg} \frac{s}{2} + \frac{1}{\cos\frac{s}{2}}\Big|\right)\Big|_x^\pi.
\end{align*}

II
\begin{align*}
&-\int \limits_{\theta_\Delta}^{\pi}\frac{\cos(s)\left(\sin \frac{s}{2} + \frac{1}{3} \sin \frac{3s}{2}\right)}{\cos\frac{s}{2}}ds = -\int \limits_{\theta_\Delta}^{\pi}\frac{\cos(s)\left(\sin \frac{s}{2} + \frac{1}{3} \left(3\sin \frac{s}{2} - 4\sin^3 \frac{s}{2} \right)\right)}{\cos\frac{s}{2}}ds = \\
&-\int \limits_{\theta_\Delta}^{\pi}\frac{\cos(s)\left(2\sin \frac{s}{2} - \frac{4}{3}\sin^3 \frac{s}{2}\right)}{\cos\frac{s}{2}}ds = -2\int \limits_{\theta_\Delta}^{\pi}\frac{\cos(s)\sin \frac{s}{2}\cos \frac{s}{2}\left(1 - \frac{2}{3}\sin^2 \frac{s}{2}\right)}{\cos^2\frac{s}{2}}ds = \\
&-2\int \limits_{\theta_\Delta}^{\pi}\frac{\cos(s)\sin s\left(1 - \frac{2}{3}\sin^2 \frac{s}{2}\right)}{\cos s + 1}ds = -\frac{2}{3}\int \limits_{\theta_\Delta}^{\pi}\frac{\cos(s)\sin s\left(2 + \cos s\right)}{\cos s + 1}ds =
\end{align*}

\begin{align*}
& u = \cos s; \hskip0.5cm du = -\sin (s) ds
\end{align*}

\begin{align*}
& = \frac{2}{3}\int \limits_{\cos(\theta_\Delta)}^{-1}\frac{u\left(2 + u\right)}{u + 1}du = \frac{2}{3} \int \limits_{\cos(\theta_\Delta)}^{\cos\pi}\left(u + 1 - \frac{1}{u+1}\right)du = \\
& = \frac{2}{3}\left(\frac{1}{2}\cos^2 s + \cos s - \ln\Big|\cos s +1\Big|\right) \Big|_x^\pi = \\
&=\frac{2}{3}\left(\frac{1}{2}\cos^2 s + \cos s - 2\ln\Big|\sqrt{2}\cos \frac{s}{2}\Big|\right) \Big|_x^\pi = \\
&=\frac{2}{3}\left(\frac{1}{2}\cos^2 s + \cos s - 2\ln\Big|\cos \frac{s}{2}\Big|\right) \Big|_x^\pi.
\end{align*}

I+II
\begin{align*}
&\left(\frac{8}{3} \sin(\frac{s}{2}) - \frac{4}{3}\ln \Big|\operatorname{tg} \frac{s}{2} + \frac{1}{\cos\frac{s}{2}}\Big|\right)\Big|_x^\pi + \frac{2}{3}\left(\frac{1}{2}\cos^2 s + \cos s - 2\ln\Big|\cos \frac{s}{2}\Big|\right) \Big|_x^\pi = \\
&\frac{1}{3}\left(8 \sin(\frac{s}{2}) - 4\ln \Big|\frac{\sin \frac{s}{2} + 1}{\cos\frac{s}{2}}\Big| + \cos^2 s + 2\cos s - 4\ln\Big|\cos \frac{s}{2}\Big|\right) \Big|_x^\pi = \\
&\frac{1}{3}\left(8 \sin(\frac{s}{2}) - 4\ln \Big|\sin \frac{s}{2} + 1\Big| + \frac{1}{2}\cos 2s + \frac{1}{2}  + 2\cos s\right) \Big|_x^\pi = \\
& = \frac{6\frac{1}{2} - 4\ln 2}{3} - \frac{1}{3}\left(8 \sin(\frac{\theta_\Delta}{2}) - 4\ln \Big|\sin \frac{\theta_\Delta}{2} + 1\Big| + \frac{1}{2}\cos 2\theta_\Delta + 2\cos \theta_\Delta\right).
\end{align*}
Hence,
\begin{align}
&S_{2}(\theta_\Delta) = \frac{K_0^2(6\frac{1}{2} - 4\ln 2) - K_0^2\left(8 \sin(\frac{\theta_\Delta}{2}) - 4\ln \Big|\sin \frac{\theta_\Delta}{2} + 1\Big| + \frac{1}{2}\cos 2\theta_\Delta + 2\cos \theta_\Delta\right)}{6\sqrt{{b} K_0}\cos\frac{\theta_\Delta}{2}} - \label{rel:approxEps2} \\
& - \frac{K_0^2 \left(\frac{2}{3} - \sin \frac{\theta_\Delta}{2} - \frac{1}{3} \sin \frac{3\theta_\Delta}{2}\right)^2}{4\sqrt{{b} K_0} \cos^3 \frac{\theta_\Delta}{2}}. \nonumber
\end{align}
In $\theta_\Delta=0$
\begin{equation}
S_{2}(0) = \frac{2 K_0^2(1 - \ln 2)}{3 \sqrt{{b} K_0}} - \frac{1 K_0^2}{9 \sqrt{{b} K_0}} = \frac{K_0^2(5 - 6\ln 2)}{9 \sqrt{{b} K_0}}. \label{rel0:approxEps2}
\end{equation}

Hence, $S_0(\theta_\Delta)$, $S_1(\theta_\Delta)$, $S_2(\theta_\Delta)$ are evaluated (equations (\ref{rel:approxEps0}), (\ref{rel:approxEps1}) and (\ref{rel:approxEps2}), correspondingly). I. e. the first and the second approximations $\hat{S}_1(\theta_\Delta, a)$, $\hat{S}_2(\theta_\Delta, a)$ of separatrix $S(\theta_\Delta, a)$ are found.
Furthermore, using (\ref{rel0:approxEps0}), (\ref{rel0:approxEps1}) and (\ref{rel0:approxEps2}) the following relations are valid:
\begin{align*}
&\hat{S}_1(0, a) = 2 \sqrt{b K_0} + a \frac{2K_0}{3}, \\
&\hat{S}_2(0, a) = \hat{S}_2(0) = 2 \sqrt{b K_0} + a \frac{2K_0}{3} + a^2 \frac{K_0(5 - 6\ln 2)}{9 b}\sqrt{K_0 b}.
\end{align*}

\section*{\uppercase{Acknowledgements}}
 This work was supported by the Russian Scientific Foundation 
 and Saint-Petersburg State University. 
 The authors would like to thank Roland~E.~Best,
 the founder of the Best Engineering Company, Oberwil, Switzerland
 and the author of the bestseller on PLL-based circuits \cite{Best-2007}
 for valuable discussion. 



\newcommand{\noopsort}[1]{} \newcommand{\printfirst}[2]{#1}
  \newcommand{\singleletter}[1]{#1} \newcommand{\switchargs}[2]{#2#1}

\end{document}